\documentclass[reqno]{amsart}
\usepackage{amsmath, amsthm, graphicx, amssymb,verbatim}
\theoremstyle{plain}
\newtheorem{theorem}{Theorem}

\newtheorem{lemma}{Lemma}
\usepackage{setspace}
\usepackage{enumitem}

\newcommand\TT{\rule{0pt}{5ex}} 

\begin{document}
\bibliographystyle{amsplain}
\title{The Farey Sieve}
\author[S. Guthery]{Scott B. Guthery}
\email{sbg@acw.com}

\begin{abstract}
An elementary method for computing various prime sequences using the sequence of Farey sequences is described.
\end{abstract}
\parindent 10pt
\maketitle
\section{Introduction}
The Farey sieve associates counts with fractions in Farey sequences and uses these counts to characterize
prime numbers.\footnote{\TT ``It is one of the clich\'es of mathematical history that it is disastrously easy even for an amateur to infer by induction results in the theory of numbers whose proof or disproof may be a century or more beyond the current capacities of mathematics.''E.T. Bell \cite[p. 178]{Bell:1947}}

\subsection{A Recursion for the Farey Sequence}
The Farey sequence $F_m$ is the sequence of irreducible rational numbers $\{\frac{n}{d}\}$ with $0 \leq n \leq d \leq m$ arranged in increasing order.

If a rational number $f$ is represented as an ordered pair $(n,d)$ where $n$ is the numerator and $d$ is the denominator of $f$, then the Farey sequence of order $m$ can be computed using the following recursion due to d'Ocagne \cite{Ocagne:1885a}:
\begin{equation}
f_{i+1}=\left[\frac{d_{i-1}+m}{d_{i}}\right]f_{i}-f_{i-1}\notag
\end{equation}
\noindent where $f_{1} =(0,1)$ and $f_{2}=(1,m)$ .

\subsection{A Recursion for the Sequence of Farey Sequences}
The original method for building the Farey sequence $F_{m+1}$ from $F_m$, the one described by Charles Haros (\cite{Haros:1801}, \cite{Haros:1802}) fifteen years before it was noticed by its namesake, John Farey \cite{Farey:1816}, is to insert $\frac{a+c}{b+d}$ between $\frac{a}{b}$ and $\frac{c}{d}$ whenever $b+d = m+1$. $\frac{a+c}{b+d}$ is the mediant of $\frac{a}{b}$~and~$\frac{c}{d}$.

While this is a handy method of generating Farey sequences using paper and pen -- which is what Charles Haros was doing when he discovered this application of the mediant -- the method does not lend itself to the analysis of the properties of Farey sequences as a function of the order, $m$.

Another way to generate $F_{m+1}$ directly from $F_m$ is keep track of where and when new elements appear. Let $f=\frac{n_f}{d_f}$ be an element of $F_m$ and define
\begin{equation}\label{SIGMA}
\sigma_f(m) = d_f-\left((m-n_f)\bmod d_f\right)\text{.}
\end{equation}

$\sigma_f(m)$ is a periodic counter associated with the fraction $f$.  It starts at $d_f-n_f$ when $f$ first appears in the sequence of Farey sequences, namely in $F_{d_f}$.  As $m$ increases it counts down to $1$ and then starts over at $d_f-1$. Whenever the counter resets, that is when it changes from $1$ in $F_{m^*}$ to $d_f-1$ in $F_{m^*+1}$ a new fraction appears on the right of $f$ in $F_{m^*+1}$, in particular one with denominator $m^*+1$.  Otherwise, the next larger fraction is the same.

\noindent Using $\sigma_f$ the Farey set $F_m$ can be expressed recursively in terms of $F_{m-1}$ as follows:
$$
F_m = F_{m-1} \bigcup \{(m-d_f,m) \mid f\in F_{m-1} \text{\, and \,} \sigma_f(m-1)=1\}
$$
\noindent where $F_1 = \{(0,1)\}$.  This isn't a closed form for $F_{m+1}$ but it offers a slight improvement on Haros' step-at-a-time mediant algorithm.  In particular based on the current sequence it lets us look ahead a number of steps to say when and where new fractions will be arriving.

The following properties of $\sigma_f$ follow directly from its definition:
\begin{enumerate}[label=\roman*)]
\item \textit{Starting Value:} $\sigma_{n_f/d_f}(d_f) = d_f - n_f$
\item \textit{Unique Values:} $\{\sigma_f(m) | d_f = m \} = \{ n_f | d_f = m\}$
\item \textit{Cycle Time:} If $\sigma_f(m)=c$, then $\sigma_f(m + k d_f)=c$ for $k \geq [1-m/d_f]$
\item \textit{Reset:} If $\sigma_{f}(m-1)=1$, then $\sigma_{f}(m)=d_f$.
\item \textit{Steps to Reset:} If $\sigma_f(m)=k$, then $\sigma_f(m+k-1)=1$.
\end{enumerate}

\noindent The two most useful properties in the following are properties $2$ and $5$. Property $2$ says that the counter values associated with the fractions having the same denominator are all different.  Property $5$ says that $\sigma_f(m)$ counts the number of steps forward from $m$ until a new (greater) fraction appears to the right of $f$.

We will see that if moving from $m$ to $m+1$ causes the reset of a fraction of every denominator from $2$ to $m$ then $m+1$ is a prime number.

\subsection{Functions of the Sequence of Farey Sequences}
The above recursion for generating the sequence of Farey sequences may be used to describe properties of Farey sequences as a function of the order of the sequence.

For example, the distance between a Farey fraction $f$ with denominator $d$ in $F_m$ and the next larger Farey fraction in $F_m$ as a function of $m$ is given by
$$
\left(\left[\frac{m-\sigma_f(m)}{d}\right]d^2 + \sigma_f(m)\, d\right)^{-1}\text{.}
$$

As another example, one can express the order index of an element $f$ of the Farey sequence $F_m$ as a function of $F_m$ and $m$ as follows. Let $i_f$ denote the order index of $f$ in $F_{d_f}$.  Then the order index of $f$ in $F_m$ for $m>d_f$ is given by
$$
I_m(f) = i_f+\sum_{g<f}\left[\frac{m-\sigma_g}{d_g}\right]- \left[\frac{d_f-\sigma_g}{d_g}\right]
$$

The result of Franel \cite{Franel:1924} regarding the equivalence of a property of the sequence of Farey sequences and the Riemann Hypothesis can be written using the above notation as
$$
\sum_{f \in F_m}\left(\frac{I_m(f)}{\Phi(m)+1}-f\right)^2 = \mathcal{O}(m^{-1+\epsilon})
$$
where $\Phi$ is the Euler summatory function.

Franel's result says intuitively that the Farey sequence encodes the prime numbers and thus it is not surprising that the Farey sequence can be used to study properties of the prime numbers.  What is perhaps a bit surprising is that the Farey sequence encodes the prime numbers without ever mentioning them.

\subsection{The Farey Sieve}
For a Farey fraction $f$, let
\begin{equation}
m_f(c,k) =
\begin{cases}
k d_f+n_f-c & c \leq n_f \\
(k+1) d_f + n_f - (1+ (c-1 \bmod d)&c>n_f
\end{cases}\notag
\end{equation}
$m_f(c,k)$ is the $k^{th}$ value of $m$ starting at $d_f$ for which $\sigma_f(m) = c$. Define
$$
M_c(d) = \bigcup_{d_f = d}\{m_f(c,k)\}_{k=0}^{\infty}\text{.}
$$
\noindent $M_c(d)$ is the (countably infinite) set of all $m$ such that there is a Farey fraction $f$ with denominator $d$ and $\sigma_f(m) = c$ in $F_m$.
Finally, define
$$
\mathcal M_c(\hat{d}) = \bigcap_{d=1}^{\hat{d}} M_c(d)\text{.}
$$
$\mathcal M_c(\hat{d})$ is the (countably infinite) set of all $m$ such that for each $d$, $1 \leq d \leq \hat{d}$, there is a Farey fraction $f$ with denominator $k$ and $\sigma_f(m)= c$ in $F_m$.

In these collections of $m$ values which particular Farey fraction with denominator $d$ satisfies a condition is immaterial.  All that matters is that there is one.   From property $2$ above there is at most one.

\section{Prime Generation}
But before we can compute $M$ and $ \mathcal M$ with only primes we must show that the sequence of primes can be constructed using the sequence of Farey sequences. The only {\it a priori} property of the primes employed in this construction is a weak one, Bertrand's Postulate (\cite{Bertrand:1845}, \cite{Chebyshev:1850}, \cite{Ramanujan:1919}, \cite{Erdos:1934}).  This postulate states that there is always a prime between $n$ and $2n$.  In particular, the sequence of prime numbers is generated without appealing to the notion of factorization.

\begin{lemma}\label{GUTHERY1}
$p+1$ is a prime if and only if $p \in \mathcal M_1(p)$ .
\end{lemma}
\begin{proof}
$p \in \mathcal M_1(p)$ if and only if there is a fraction $f$ with denominator $d$ and $\sigma_f(p)=1$ in $F_{p}$ for
each $d$, $1 \leq d \leq p$. In this case $|F_{p+1}|-|F_p| = p = \varphi(p+1)$ and  $\varphi(p+1)=p$ if and only if $p+1$ is a prime.
\end{proof}

\begin{theorem}\label{GUTHERY2}
Let $p$ be a prime and set
$$
{\bar p}= min\, \mathcal M_1(p)\text{.}
$$
Then ${\bar p}+1$ is the next prime after $p$.
\end{theorem}
\begin{proof}
Consider $d \in \mathcal M_1(p)$. From the definition of $\mathcal M_1(p)$, $1 \leq d \leq p$ and there is an $f$ with denominator $d$ in $F_{\bar p}$ such that $\sigma_f({\bar p}) = 1$.  If in addition $d \leq [{\bar p}/2]$ then $\sigma_f({\bar p}-d) = 1$.

If $m = {\bar p}-d+1$ then there is $f' \in F_m$ with $\sigma_{f'}(m) = d$. Then $\sigma_{f'}(m+d-1) = \sigma_{f'}({\bar p} = 1$.

To summarize, for each $m$ such that ${\bar p}/2+1 < m \leq {\bar p}$, there is an $f$ with denominator $m$ in $F_{\bar p}$ such that $\sigma_f({\bar p})=1$.

Since ${\bar p}/2 < p$ by Bertrand's Postulate  we have that there is a $f$ in $F_{\bar p}$ with $\sigma_f({\bar p})=1$ and denominator $d$ every $d$ between $1$ and ${\bar p}$. Therefore ${\bar p} \in \mathcal M_1({\bar p})$ and ${\bar p}+1$ is a prime by Lemma \ref{GUTHERY1}.

Since ${\bar p} = min\, \mathcal M_1(p)$, for each $m$, $p \leq m < {\bar p}$, there must be at least one $d$, $1 \leq d \leq m$, for which there is no fraction with denominator $d$ and $\sigma_f(m)=1$ in $F_m$.  Therefore, $m$ is not prime and ${\bar p}+1$ is the next prime after $p$.
\end{proof}

\section{Farey Sieving}
The following lemmas help construct efficient sieves.

\begin{lemma}
If $p$ is a prime factor of $d$ then $M_c(d) \subset M_c(p)$.
\end{lemma}
\begin{proof}
Each $m$ in $M_c(d)$ is of the form
$$
m = k_d d_f + n_f - c = k_d (p q) + n_f - c
$$
where $q=\frac{d}{p}$ and $f \in F_d$. Since $n_f \bmod p \neq 0$, $n_f = k'p + n_p$ where $1 \leq n_p < p$. As $p$ is prime $\frac{n_p}{p} \in F_p$ so
$$
m = (k_d q + k')p + n_p
$$
is in $M_c(p)$.
\end{proof}

\begin{lemma}
If $\{p_i\}$ is the set of prime factors of $d$ then
$$
M_c(d) = \left(\bigcap M_c(p_i)\right)/\{1,\ldots,d\}\text{.}
$$
\end{lemma}
\begin{proof}
From Lemma 1, $M_c(d)  \subset \left(\bigcap M_c(p_i)\right).$ If $m-c = n_i \bmod p_i$ then by the Chinese Remainder Theorem, there is a $y$ such that $m-c = y \bmod d$ so $m \in M_c(d)$.
\end{proof}

\begin{lemma}
For prime $p$, $M_c(p) =  \{p, p+1, \ldots\} / \lbrace k p - (c \bmod p) \rbrace_{k=1}^\infty$
\end{lemma}
\begin{proof}
There are $p-1$ fractions with denominator $p$ each with a different $\sigma_f(m)$ value that is between $1$ and $p$ therefore exactly one of these values is always missing; namely $p - (m \bmod p)$.
\end{proof}

As a result of these three lemmas we have
$$
\mathcal M_c(d) = \bigcap_{p \leq d} M_c(p) = I \hspace{2mm} / \bigcap_{p\; \leq \; d}\lbrace k p - (c \bmod p) \rbrace_{k=1}^\infty\text{.}
$$

\section{Twin Prime Generation}
The twin prime conjecture is a special case of a conjecture due to Polignac \cite{Polignac:1849c} who posited that there are an infinite number prime pairs, $p_1$ and $p_2$, such that $p_2 = p_1 + k$ where $k$ is a fixed even number. The case when $k = 2$ is called the twin prime conjecture.
\begin{lemma}\label{GUTHERY3}
$p+1$ is the lesser of a twin prime if and only if
$$
p \in \mathcal M_1(p) \bigcap \mathcal M_3(p)\text{.}
$$
\end{lemma}
\begin{proof}
\noindent $p$ is a prime since $p-1 \in \mathcal M_1(p-1)$.

\noindent Further since $([p/2],p,2) \in F_p$ we have $([p/2],p,1) \in F_{p+1}$ along with $(1,1)$ and $(p+1,1)$. The condition
$$
p-1 \in \mathcal M_3(p-1)
$$
ensures that $(d, 1) \in F_{p+1}$ for $2 \leq d \leq p-1$.  Therefore $(d, 1) \in F_{p+1}$ for $1 \leq d \leq p+1$. That is,
$$
p+1 \in \mathcal M_1(p+1)
$$
so $p+2$ is a prime.
\end{proof}

\textbf{Twin Prime Conjecture}
\textit{Let $p$ be the lesser of a twin prime and set
$$
{\bar p}= \min \left( \mathcal M_1(p) \bigcap \mathcal M_3(p)\right)\text{.}
$$
Then ${\bar p}+1$ is the lesser of the next twin prime after $p$.}

From the above lemmas, for prime $p$ and $q$,
$$
\mathcal M_1(p) \bigcap \mathcal M_3(p) =
I\hspace{2mm} /  \bigcap_{\text{prime } q \leq p}\left[\lbrace k q - 1) \rbrace_{k=1}^\infty \bigcup \lbrace k q - 3) \rbrace_{k=1}^\infty\right]
$$
so $M_1(p) \bigcap \mathcal M_3(p)$ is never null and ${\bar p}$ always exists.

\section{Prime Constellations}
A sequence of integers $a_1 < a_2 < ... < a_n$ is called an {\it (admissible) prime constellation} if for each prime $p$, there is some residue class modulo $p$ which contains none of the $a_i$. Hardy and Littlewood \cite{Hardy:1922} conjectured that if $\lbrace a_i \rbrace_{i=1}^k$ is a prime constellation then there are infinitely many integers $n$ such that $n+a_1, \ldots, n+a_k$ are all prime.

Table~\ref{TABLE1} lists some examples of prime constellations. The first entry in this table describes the prime numbers and the second entry describes prime pairs. Table~\ref{TABLE2} lists some instances of the fourth entry in Table \ref{TABLE1}, the prime constellation $(0, 2, 6, 8)$. The sequence of primes
$$
(284723,284729,284731,284737,284741,284743,284747,284749)
$$
is described by the eighth entry, $\{0, 6, 8, 14, 18, 20, 24, 26\}$. For many more examples see \cite{Erdos:1988}, \cite{Guy:2004a}, \cite{Odlyzko:1999}, and  \cite{Schinzel:1958a} as well as various web sites devoted to prime numbers.

\begin{table}[!ht]
\centering
\begin{tabular}{c}
$\{0 \}$\\
$\{0, 2 \}$\\
$\{0, 4, 6\}$\\
$\{0, 2, 6, 8 \}$\\
$\{0, 4, 6, 10, 12)\}$\\
$\{0, 4, 6, 10, 12, 16\}$\\
$\{0, 2, 8, 12, 14, 18, 20\}$\\
$\{0, 6, 8, 14, 18, 20, 24, 26\}$\\
$\{0, 2, 6, 12, 14, 20, 24, 26, 30\}$\\
$\{0, 2, 6, 12, 14, 20, 24, 26, 30, 32\}$
\end{tabular}
\vspace{5mm}\caption{Examples of Prime Constellations}\label{TABLE1}
\end{table}
\begin{table}[!ht]
\centering
\begin{tabular}{c}
$\{5,7,11,13\}$\\
$\{11,13,17,19\}$\\
$\{821,823,827,829\}$\\
$\{3251,3253,3257,3259)\}$\\
\end{tabular}
\vspace{5mm}\caption{Examples of the Prime Constellation $(0, 2, 6, 8)$}\label{TABLE2}
\end{table}

In this note we use the Farey sequence to easily compute instances of a given prime constellation. The method is then used to demonstrate that for any prime constellation there are an infinite number of such instances.

The sequence of $k$-tuples satisfying a prime constellation can be generated by taking the $\mathcal M$ intersection over the constellation.  The above result for twin primes can be easily extended to include prime constellations.

\textbf{Prime Constellation Conjecture}
Let $p$ be the minimum prime of a $k$-tuple described by the prime constellation $\lbrace a_i\rbrace_{i=1}^n$ and set
$$
{\bar p}= \min \left(\bigcap_{i=1}^n \mathcal M_{a_i+1}(p)\right)\text{.}
$$
Then ${\bar p}+1$ is the minimum prime in the next $k$-tuple described by the constellation

\section{Prime Function Values for a Common Argument}

Dickson conjectured in $1904$ \cite{Dickson:1904a} that given a family of linear functions with integer coefficients $a_i > 1$ and $b_i$:
$$
a_1 n + b_1, a_2 n + b_2, \ldots, a_k n + b_k
$$
then there are infinitely many integers $n > 0$ for which these are simultaneously prime unless there is a prime $p$ which divides the product of these for all $n$.

\noindent\textbf{Dickson Conjecture.}
Let $p$ be the minimum prime of a $k$-tuple of primes described by the linear coefficients $\lbrace a_i i + b_i\rbrace_{i=1}^n$ and set
$$
{\bar p}= \min \left(\bigcap_{i=1}^n \mathcal M_{1}(a_i i + b_i)\right)\text{.}
$$
Then ${\bar p}+1$ is the minimum prime in the next $k$-tuple described by the linear coefficients.

Schinzel-Sierpinski \cite{Schinzel:1958} and Bateman-Horn \cite{Bateman:1962} proposed the following generalization of the conjectures of Bunyakovskii, Goldbach, Hardy-Littlewood, Dickson, Shanks and others prime function values for a common argument.

\textbf{Hypothesis H.} Let $k$ be a positive integer and let $f_1(x), f_2(x),\ldots, f_k(x)$ be irreducible polynomials with integral coefficients and positive leading coefficients such that there is not a prime $p$ which divides the product
$$
f_1(m)\cdot f_2(m)\cdot \ldots \cdot f_k(m)
$$
for every integer $m$.  Then there exists an integer $n$ such that $f_1(n), f_2(n), \ldots , f_k(n)$ are all prime numbers.

\pagebreak

\section{Mathematica Code for the Farey Sieve}
The following Mathematica$^\circledR$ notebook implements the Farey sieve for arbitrary sets of functions. Table~\ref{TABLE3} lists invocations of this notebook to generate various familiar prime sequences.
\begin{verbatim}
IntegerTable = {};

FunctionSieve[ems_, f_] :=
 Module[{i, iMax, d, l = {}, emsMax = ems[[Length[ems]]]},
  d = Map[f, IntegerTable] - 1;
  iMax = LengthWhile[d, # <= emsMax &];
  For[i = 1, i <= iMax, i++,
   If[MemberQ[ems, d[[i]]],
     l = Append[l, i];
     ];
   ];
  l
  ]

FareySieve[f_, n_: 6, size_: 1000] :=
 Module[{i, j, d = 2, l, tp = 2, ems, fun, k = n, p, t},
  l = IntegerTable = Table[i, {i, 2, size}];
  While[k > 0,
   ems = Fms[d];
   For[i = 1, i <= Length[f], i++,
    l = Intersection[l, FunctionSieve[ems, f[[i]]]];
    ];
   If[Length[l] == 0,
    Print["Null Intersection"];
    Break[]
    ];
   If[ l[[1]] + 1 != tp,
    tp = l[[1]] + 1;
    t = {};
    For[i = 1, i <= Length[f], i++,
     p = f[[i]][tp];
     t = Append[t, {p, PrimeQ[p]}];
     ];
    Print[t];
    k--;
    ];
   d++;
   ];
  ]
\end{verbatim}

\begin{table}[!ht]
\centering
\begin{tabular}{l l}
  Name & Invocation\\
  \hline
  Primes & \texttt{FareySieve[\{\#\&\}]}\\
  Twin Primes & \texttt{FareySieve[\{\#\&,\#+2\&\}]}\\
  Prime Constellation & \texttt{FareySieve[\{\#\&,\#+2\&,\#+6\&\}]}\\
  Sophie Germain Primes & \texttt{FareySieve[\{\#\&,2\#+1\&\}]}\\
  Gaussian Primes & \texttt{FareySieve[\{\#\&,4\#+3\&\}]}\\
  Cunningham Chain & \texttt{FareySieve[\{\#\&,2\#+1\&,4\#+3\&\}]}\\
  Dickson Chain & \texttt{FareySieve[\{\#\&,2\#+1\&,3\#+4\&\}]}\\
  Star Primes & \texttt{FareySieve[\{6\#(\#-1)+1\&\}]}\\
  Shanks Primes & \texttt{FareySieve[\{\#\^{}2+1\&\}]}\\
  Shanks Twins & \texttt{FareySieve[\{(\#-1)\^{}2+1 \&,(\#+1)\^{}2+1 \&\}]}\\
  Hardy-Littlewood Primes & \texttt{FareySieve[\{\#\^{}2+\#+1\&\}]}\\
  Thabit Primes & \texttt{FareySieve[\{3*2{}\^{}\#-1\&\}]}\\
  Wagstaff Primes & \texttt{FareySieve[\{(2\^{}\#+1)/3\&\}]}\\
  Safe Primes & \texttt{FareySieve[\{\#\&,(\#-1)/2\&\}]}\\
  Proth Primes & \texttt{FareySieve[\{2\^{}\#+1\&\}]}\\
  Kynea Primes & \texttt{FareySieve[\{(2\^{}\#+1)\^{}2-2\&\}]}\\
  Mersenne Primes & \texttt{FareySieve[{2\^{}\#-1\&}]}\\
  Centered Heptagonal Primes & \texttt{FareySieve[\{(7\#\^{}2-7\#+2)/2\&\}]}\\
  Centered Square Primes & \texttt{FareySieve[\{\#\^{}2+(\#+1)\^{}2\&\}]}\\
  Centered Triangular Primes & \texttt{FareySieve[\{(3\#\^{}2+3\#+2)/2\&\}]}\\
  Carol Primes & \texttt{FareySieve[\{(2\^{}\#-1)\^{}2-2\&\}]}\\
  Cullen Primes & \texttt{FareySieve[\{\#2\^{}\#+1\&\}]}\\
  Double Mersenne Primes & \texttt{FareySieve[\{\#\&,2\^{}(2\^{}\#-1)-1\&\}]}\\
  Euler Primes & \texttt{FareySieve[\{2\^{}(2\^{}\#)+1\&\}]}
\end{tabular}
\vspace{5mm}\caption{Invocations of Farey Sieve for Various Prime Sequences}\label{TABLE3}
\end{table}
\pagebreak
\bibliography{PrimeRecursions}
\end{document}